\documentclass[preprint,12pt]{elsarticle}
\usepackage{amsfonts}
\usepackage{amsmath}
\usepackage{mathrsfs}
\usepackage[applemac]{inputenc}
\usepackage{amssymb}
\usepackage{times, ulem, amsfonts}
\usepackage{mathpazo, amsmath, amssymb, amsthm,color}

\newcommand{\p}{{\mathbb P}}

\def\virgp{\raise 2pt\hbox{,}}
\def\({\left(}
\def\){\right)}

\def\<{\langle}
\def\>{\rangle}

\def\({\left(}
\def\){\right)}
\theoremstyle{plain}

\numberwithin{equation}{section}
\newtheorem{theorem}{Theorem}[section]

\newtheorem{lemma}[equation]{Lemma}

\newtheorem{definition}[equation]{Definition}

\theoremstyle{remark}
\newtheorem{remark}[theorem]{Remark}

\newcommand{\R}{{\mathbb R}}
\def\div{ \hbox{\rm div}\,  }
\newcommand\Z{{\mathbb{Z}}}

\def\aa{\phi}
\def\ddj{\dot \Delta_j}

\begin{document}
\bibliographystyle{plainmma}
\baselineskip=24pt
\footnote{
Email Addresses: zhaixp@szu.edu.cn (X. Zhai).}
\begin{center}
{\Large\bf  Optimal decay for the   $n$-dimensional incompressible Oldroyd-B model without damping mechanism}\\[1ex]

Xiaoping Zhai \\[1ex]
 School  of Mathematics and Statistics, Shenzhen University,
 Shenzhen 518060, China
\end{center}

\centerline{\Large\bf Abstract}
\bigskip
\noindent
By a new energy approach involved in the high frequencies and low frequencies decomposition in the Besov spaces, we obtain the  optimal decay for the   incompressible Oldroyd-B model without damping mechanism in $\R^n$ ($n\ge 2$). More precisely, let $(u,\tau)$ be the global small solutions constructed in \cite{zhaixiaoping},
  we prove for any  $(u_0,\tau_0)\in{\dot{B}_{2,1}^{-s}}(\R^n)$  that
\begin{eqnarray*}
\big\|\Lambda^{\alpha}(u,\Lambda^{-1}\p\div\tau)\big\|_{L^{q}} \le C \(1+t\)^{-\frac n4-\frac {(\alpha+s)q-n}{2q}}, \quad\Lambda\stackrel{\mathrm{def}}{=}\sqrt{-\Delta},
\end{eqnarray*}
 with
$\frac n2-1<s<\frac np, $
$2\leq p \leq \min(4,{2n}/({n-2})),\  p\not=4\ \hbox{ if }\ n=2,$
and
$p\leq q\leq\infty$,
$\frac nq-\frac np-s<\alpha \leq\frac nq-1$.
The proof relies heavily on the special dissipative  structure of the equations and some   commutator estimates and various interpolations
between Besov type spaces.
The method also works for other parabolic-hyperbolic systems in which the   Fourier splitting technique is invalid.

\date{\today}


\setcounter{section}{1}
\setcounter{theorem}{0}
\section*{\Large\bf 1. Introduction and the main result}
The  incompressible Oldroyd-B model without damping mechanism in $\R^n$ can be written as:
\begin{eqnarray}\label{m}
\left\{\begin{aligned}
&\partial_t\tau + u\cdot \nabla \tau   + F(\tau, \nabla u)  =  D (u),\\
&\partial_t u+ u\cdot\nabla u-\Delta u+\nabla{p}= \div \tau,\\
&\div u =0,\\
&(u,\tau)|_{t=0}=(u_0,\tau_0),
\end{aligned}\right.
\end{eqnarray}
where
 $u=(u_{1},u_{2},\cdot\cdot \cdot, u_{n})$ denotes the velocity, ${p}$ is the scalar pressure of fluid. $\tau=\tau_{i,j}$ is the non-Newtonian part of stress tensor which can be seen as a symmetric matrix here. $D(u) = \frac{1}{2} \big( \nabla u + (\nabla u)^{T} \big)$ is the symmetric part of $\nabla u$
and $ F$ is a given bilinear form which can be chosen as
\begin{equation*}
F(\tau, \nabla u)= \tau \Omega(u) - \Omega(u) \tau + b(D(u) \tau + \tau D(u)),
\end{equation*}
where $b$ is a parameter in $[-1,1]$, $\Omega(u) = \frac{1}{2} \big( \nabla u - (\nabla u)^{T} \big)$ is the skew-symmetric part of $\nabla u$.

The above Oldroyd-B model presents a typical constitutive law which does not obey the Newtonian law (a linear relationship between stress and the gradient of velocity
in fluids). Such non-Newtonian property may arise from the memorability of some fluids. Formulations about viscoelastic flows of Oldroyd-B type are first introduced by Oldroyd \cite{Oldroyd} and are extensively discussed in \cite{bird}.
One can find
 the derivation of \eqref{m} in \cite{lin2012}, here we omit it.

The mathematical theory of Oldroyd-B model  is an old subject, see \cite{chemin}--\cite{GS2}, \cite{LMZ}--\cite{Oldroyd}, \cite{zhaixiaoping}--\cite{zuiruizhao}.
Here, we only recall some results about Oldroyd-B model without damping mechanism. In fact, when neglecting the damping term in the stress tensor equation, \eqref{m} reduces to be a parabolic-hyperbolic system. Duo to lack of smoothing effect of $\tau$, it's difficult to get the global solutions directly. Luckily,  by exploiting the good structure of the system, we can obtain some hidden dissipation about $\tau.$
Based on the above analysis,
by constructing the time-weighted energies,
Zhu  \cite{zhuyi} obtained the global small  solutions to \eqref{m} in $\R^3.$
This result was extended by
Chen and Hao  \cite{chenqionglei} to the  $L^2$ type Besov spaces in $\R^n.$ The first author  of the present paper in \cite{zhaixiaoping} generalized the result of
\cite{chenqionglei} to the  $L^p$ framework which the highly oscillating
initial velocity are allowed.

Denote $\p=\mathcal{I}-\mathcal{Q}:=\mathcal{I}-\nabla\Delta^{-1}\div$ and
$$f^\ell\stackrel{\mathrm{def}}{=}\dot{S}_{j_0+1}f\quad\hbox{ and }\quad f^h\stackrel{\mathrm{def}}{=}f-f^\ell$$ for some fix integer $j_{0}\ge0$,

The author in \cite{zhaixiaoping} obtained
the following theorem:
\begin{theorem}\label{th1} (see \cite{zhaixiaoping})
Let   $n\ge 2$ and
\begin{equation*}
2\leq p \leq \min(4,{2n}/({n-2}))\quad\hbox{and, additionally, }\  p\not=4\ \hbox{ if }\ n=2.
\end{equation*}
 For any  $(u_0^\ell,\tau_0^\ell)\in \dot{B}_{2,1}^{\frac n2-1}(\R^n)$, $u_0^h\in \dot{B}_{p,1}^{\frac np-1}(\R^n)$,  $\tau_0^h\in \dot{B}_{p,1}^{\frac np}(\R^n)$ with $\div u_0 = 0$.
 If  there exists a positive constant $c_0$ such that,
\begin{align*}
{\mathcal X}_{0} \stackrel{\mathrm{def}}{=}\|(u_0,\tau_0)\|^\ell_{\dot{B}_{2,1}^{\frac {n}{2}-1}}+ \|u_0\|^h_{\dot B^{\frac  np-1}_{p,1}}+\|\tau_0\|^h_{\dot B^{\frac  np}_{p,1}}\leq c_0,
\end{align*}
then
the system \eqref{m} has a unique global solution $(u,\tau)$ so that for any $T>0$
\begin{align*}
&u^\ell\in C_b([0,T );{\dot{B}}_{2,1}^{\frac {n}{2}-1}(\R^n))\cap L^{1}
([0,T];{\dot{B}}_{2,1}^{\frac n2+1}(\R^n)),\\
&
\tau^\ell\in C_b([0,T );\dot{B}_{2,1}^{\frac n2-1}(\R^n)), \quad(\Lambda^{-1}\p\div\tau)^\ell\in L^{1}
([0,T];{\dot{B}}_{2,1}^{\frac n2+1}(\R^n)),\\
&u^h\in C_b([0,T );{\dot{B}}_{p,1}^{\frac {n}{p}-1}(\R^n))\cap L^{1}
([0,T];{\dot{B}}_{p,1}^{\frac np+1}(\R^n)),\\
&
\tau^h\in C_b([0,T );\dot{B}_{p,1}^{\frac np}(\R^n)), \quad (\Lambda^{-1}\p\div\tau)^h\in L^{1}
([0,T];{\dot{B}}_{p,1}^{\frac np}(\R^n)).
\end{align*}
Moreover,  there exists some constant $C=C(p,n)$ such that
\begin{equation*}
{\mathcal X}(t)\leq C {\mathcal X}_{0},
\end{equation*}
\begin{align*}
with \quad
{\mathcal X}(t)\stackrel{\mathrm{def}}{=}&\|(u,\tau)\|^\ell_{\widetilde{L}^\infty_t(\dot{B}_{2,1}^{\frac{n}{2}-1})}
+\|u\|^h_{\widetilde{L}^\infty_t(\dot{B}_{p,1}^{\frac{n}{p}-1})}
+\|\tau\|^h_{\widetilde{L}^\infty_t(\dot{B}_{p,1}^{\frac{n}{p}})}+\|u\|^h_{L^1_t(\dot{B}_{p,1}^{\frac{n}{p}+1})}
\\&+\|(u,(\Lambda^{-1}\p\div\tau))\|^\ell_{L^1_t(\dot{B}_{2,1}^{\frac{n}{2}+1})}+\|\Lambda^{-1}\p\div\tau\|^h_{L^1_t(\dot{B}_{p,1}^{\frac{n}{p}})}.
\end{align*}
\end{theorem}
The natural next step is to look for a more accurate description of the long time behavior
of the solutions.
 As  there is no dissipation in the $\tau$ equation, the usual Fourier splitting technique can not be used here.
 The spectral analysis for the linearized system may be valid. Here, we present another new  pure energy method which motivated by \cite{guoyan}, \cite{xinzhouping} to  get the optimal decay of the solutions.
 Considering the linear system of \eqref{m}, one can find $u$ and $\p\div\tau$ satisfy the following damped wave equation:
 \begin{equation*}
W_{tt}-\Delta W_{t}-\frac12\Delta W=0.
\end{equation*}
  Thus, we only expect to get the  decay of  $u$ and the partial decay in $\tau$, namely $\p\div\tau$.

 Now, we  state the main result of the paper:
\begin{theorem}\label{th2}
Let ~$(u,\tau)$ be the global small solutions addressed by Theorem \ref{th1}.
If in addition   $(u_0,\tau_0)\in{\dot{B}_{2,1}^{-s}}(\R^n)$  with
$\frac n2-1<s<\frac np. $
For any
$p\leq q\leq\infty$ and
$\frac nq-\frac np-s<\alpha \leq\frac nq-1$, there holds
\begin{eqnarray}\label{xiaoxiao}
\big\|\Lambda^{\alpha}(u,\Lambda^{-1}\p\div\tau)\big\|_{L^{q}} \le C \(1+t\)^{-\frac n4-\frac {(\alpha+s)q-n}{2q}}.
\end{eqnarray}
\end{theorem}
\begin{remark}\label{yaya3}
Let $p=q=2$, one can deduce from \eqref{xiaoxiao} that
\begin{eqnarray*}
\big\|\Lambda^{\alpha}(u,\Lambda^{-1}\p\div\tau)\big\|_{L^{2}} \le C \(1+t\)^{-\frac s2-\frac {\alpha}{2}},
\end{eqnarray*}
which coincides with the  heat flows, thus our decay rate is optimal in some sense.
\end{remark}

\setcounter{equation}{0}
\setcounter{section}{2}
\section*{\large\bf 2. Preliminaries}
For  readers' convenience, in this section, we list some basic knowledge on  Littlewood-Paley  theory.

\begin{definition}
 Let us consider a smooth function~$\varphi $
on~$\R,$ the support of which is included in~$[\frac34,\frac83]$ such that
$$
\forall
 \tau>0\,,\ \sum_{j\in\Z}\varphi(2^{-j}\tau)=1, \quad\hbox{and}\quad \chi(\tau)\stackrel{\mathrm{def}}{=} 1 - \sum_{j\geq
0}\varphi(2^{-j}\tau) \in {\mathcal D}([0,4/3]).
$$
Let us define
$$
\ddj u={\mathcal F}^{-1}(\varphi(2^{-j}|\xi|)\widehat{u}),
 \quad\hbox{and}\quad \dot{S}_ju={\mathcal F}^{-1}(\chi(2^{-j}|\xi|)\widehat{u}).
$$
Let $p$ be in~$[1,+\infty]$ and~$s$ in~$\R$, $u\in\mathcal{S}'(\R^n)$. We define the Besov norm by
$$
\|u\|_{{\dot{B}^s_{p,1}}}\stackrel{\mathrm{def}}{=}\big\|\big(2^{js}\|\ddj
u\|_{L^{p}}\big)_j\bigr\|_{\ell ^{1}({\mathop{\mathbb Z\kern 0pt}\nolimits})}.
$$
We then define the spaces
$\dot{B}_{p,1}^s\stackrel{\mathrm{def}}{=}\{u\in\mathcal{S}'_h(\R^n),\
\|u\|_{\dot{B}_{p,1}^s}<\infty\}$, where $u\in \mathcal{S}'_h(\R^n)$ means that $u\in \mathcal{S}'(\R^n)$ and $\lim_{j\to -\infty}\|\dot{S}_ju\|_{L^\infty}=0$ (see Definition 1.26 of \cite{bcd}).
\end{definition}

In this paper, it will be suitable to split tempered distributions $u$  into  low and high frequencies.
For a fix  integer $j_{0}$ (the value of which will
follow from the proof of the main theorem), we denote
$$\left\| u\right\| _{\dot{B}_{p,1}^{s}}^{\ell} \stackrel{\mathrm{def}}{=} \sum_{j\leq
j_{0}}2^{js}\| \dot{\Delta}_{j}u\|_{L^{p}} \ \ \mbox{and} \ \ \|u\|_{\dot{B}_{p,1}^{s}}^{h}\stackrel{\mathrm{def}}{=} \sum_{j\geq j_{0}-1}2^{js}\| \dot{\Delta}_{j}u\| _{L^{p}}.$$

Let us now state some classical
properties for the Besov spaces.
\begin{lemma}\label{qianru}

\begin{itemize}
\item Let  $1\le p\le \infty$ and $s_1,\ s_2\in \R$  with $s_1>s_2$,    for any $u\in\dot{B}^{s_1}_{p,1}\cap\dot{B}^{s_2}_{p,1}(\R^n)$, there holds
\begin{align*}
\|u^\ell\|_{\dot{B}^{s_1}_{p,1}}\le C\| u^\ell\|_{\dot{B}^{s_2}_{p,1}}, \quad \|u^h\|_{\dot{B}^{s_2}_{p,1}}\le C\| u^h\|_{\dot{B}^{s_1}_{p,1}}.
\end{align*}

\item If $s_1\neq s_2$ and $\theta\in(0,1)$,
$\left[\dot{B}_{p,1}^{s_1},\dot{B}_{p,1}^{s_2}\right]_{\theta}=\dot{B}_{p,1}^{\theta
s_1+(1-\theta)s_2}$.

\item For any smooth homogeneous of degree $m\in\Z$ function $A$ on $\R^n\backslash\{0\}$, the operator $A(D)$ maps $\dot{B}^{s}_{p,1}$ in $\dot{B}^{s-m}_{p,1}$.
\end{itemize}
\end{lemma}

We are going to define the space of Chemin-Lerner (see \cite{bcd}) in which we will work, which is a
refinement of the space ${L^\lambda_{T}(\dot{B}_{p,1}^s(\mathbb{R}^n))}$.

\begin{definition}
Let $(\lambda,p) \in [1,+\infty]^2$ and $T \in (0,+\infty]$. We define
${\widetilde{L}^\lambda_{T}(\dot{B}_{p,1}^s(\mathbb{R}^n))}$ as the completion of $C([0,T ];\mathscr{S}(\mathbb{R}^n))$ by the norm
$$
\|f\|_{\widetilde{L}^\lambda_{T}(\dot{B}_{p,1}^s)}=\sum_{j\in \mathbb{Z}}2^{js}
\left(\int_0^T\|\dot{\Delta}_jf(t)\|_{L^p}^\lambda dt\right)^{\frac{1}{\lambda}}<\infty.
$$
\end{definition}

The following product estimates in Besov spaces play a key role in our analysis of the bilinear terms (see \cite{xinzhouping}).
\begin{lemma}\label{daishu}
  Let $1\leq p, q\leq \infty$, $s_1\leq \frac{n}{q}$, $s_2\leq n\min\{\frac1p,\frac1q\}$ and $s_1+s_2>n\max\{0,\frac1p +\frac1q -1\}$. For $\forall (u,v)\in\dot{B}_{q,1}^{s_1}({\mathbb R} ^n)\times\dot{B}_{p,1}^{s_2}({\mathbb R} ^n)$, we have
\begin{align*}
\|uv\|_{\dot{B}_{p,1}^{s_1+s_2 -\frac{n}{q}}}\le C \|u\|_{\dot{B}_{q,1}^{s_1}}\|v\|_{\dot{B}_{p,1}^{s_2}}.
\end{align*}

\end{lemma}

Finally, we recall the following commutator's estimate:
\begin{lemma}(Lemma 2.100 from Bahouri et al. (2011))\label{jiaohuanzi}
Let  $\nabla u\in \dot{B}_{p,1}^{\frac {n}{p}}(\R^n)$ and $v \in \dot{B}_{q,1}^{s}(\R^n)$  with  $1\leq p, q,r\leq \infty$.
For any
$$-1-n\min\left\{\frac 1p\, ,1-\frac {1}{q}\right\}<s\le\frac np,\ \ \ \ \quad \hbox{ if}\quad \div u=0,$$
 there holds
$$
\big\|\big(2^{js}\|[\dot{\Delta}_j, u\cdot \nabla ]v\|_{L^{q}}\big)_j\bigr\|_{\ell ^{1}({\mathop{\mathbb Z\kern 0pt}\nolimits})}
\le C\|\nabla u\|_{\dot{B}_{p,1}^{\frac {n}{p}}}\|v \|_{\dot{B}_{q,1}^{s}}.
$$
\end{lemma}

\setcounter{equation}{0}
\setcounter{section}{3}
\section*{\large\bf 3. Proof of the main theorem}
In this section, we prove the main Theorem \ref{th2} by a pure energy method which
is originated from the idea as in \cite{guoyan}, \cite{xinzhouping}.

Applying project operator $\p$ on both hand side of the first two equation in \eqref{m} gives
\begin{eqnarray}\label{G1}
\left\{\begin{aligned}
&\partial_t  u +\p (u\cdot \nabla u)-\Delta  u-\p \div \tau=0,\\
&\partial_t \p \div \tau+\p \div(u\cdot \nabla \tau) -\Delta u+\p \div(F(\tau, \nabla u))=0.
\end{aligned}\right.
\end{eqnarray}
Denote
$$\aa\stackrel{\mathrm{def}}{=}\Lambda^{-1}\p\div\tau,\quad  \hbox{with}\quad\Lambda\stackrel{\mathrm{def}}{=}\sqrt{-\Delta}.$$

A simple computation from \eqref{G1}  implies
\begin{eqnarray}\label{G5}
\left\{\begin{aligned}
&\partial_t  u + u\cdot \nabla u-\Delta  u-\Lambda \aa=-[\p ,u\cdot \nabla] u,\\
&\partial_t \aa+ u\cdot \nabla\aa+\Lambda u=- [\Lambda^{-1}\p\div,u\cdot \nabla ]\tau-\Lambda^{-1}\p \div(F(\tau, \nabla u)) .
\end{aligned}\right.
\end{eqnarray}
Now, we can follow the proof of  Section 3 in \cite{zhaixiaoping} or Lemma 4.1 and Lemma  4.2 in \cite{xinzhouping}
to get (we omit the details)
\begin{align}\label{A1}
&\frac{d}{dt}\big(\|(u,\aa)\|^\ell_{\dot{B}_{2,1}^{\frac n2-1}}+\|u\|^h_{\dot B^{\frac  np-1}_{p,1}}+\|\aa\|^h_{\dot B^{\frac  np}_{p,1}}\big)
+\|( u,\aa)\|^\ell_{\dot{B}_{2,1}^{\frac n2+1}}
+ \|u\|^h_{\dot B^{\frac  np+1}_{p,1}}+ \|\aa\|^h_{\dot B^{\frac  np}_{p,1}}
\nonumber\\
&\quad
\le C\big(\|(u,\tau)\|^\ell_{\dot B^{\frac  n2-1}_{2,1}}+\|u\|^h_{\dot B^{\frac  np-1}_{p,1}}+\|\tau\|^h_{\dot{B}_{p,1}^{\frac np}}\big)\big(\| (u,\tau)\|^\ell_{\dot B^{\frac  n2+1}_{2,1}}  +\|u\|^h_{\dot B^{\frac  np+1}_{p,1}}+\|\aa\|^h_{\dot B^{\frac  np}_{p,1}}\big).
\end{align}
The  following fact can be guaranteed  by Theorem \ref{th1}:
\begin{align}\label{xiaotiao}
\|(u,\tau)\|^\ell_{\dot B^{\frac  n2-1}_{2,1}}+\|u\|^h_{\dot B^{\frac  np-1}_{p,1}}+\|\tau\|^h_{\dot{B}_{p,1}^{\frac np}}\le {\mathcal X}(t)\le\mathcal{X}_{0}\ll 1
\quad \hbox{for all } t\geq0.
\end{align}

Thus absorbing all the terms in the right  to  left in \eqref{A1} gives
\begin{align}\label{A2}
&\frac{d}{dt}\big(\|(u,\aa)\|^\ell_{\dot{B}_{2,1}^{\frac n2-1}}+\|u\|^h_{\dot B^{\frac  np-1}_{p,1}}+\|\aa\|^h_{\dot B^{\frac  np}_{p,1}}\big)
+\frac12\big(\|( u,\aa)\|^\ell_{\dot{B}_{2,1}^{\frac n2+1}}
+ \|u\|^h_{\dot B^{\frac  np+1}_{p,1}}+ \|\aa\|^h_{\dot B^{\frac  np}_{p,1}}\big)\le 0.
\end{align}

Next, we want to use the interpolation inequality to get the Lyapunov-type inequality for the above energy norms.

According to \eqref{xiaotiao} and  Lemma \ref{qianru},
it's obvious  for any $\beta>1$ that
\begin{align}\label{A812}
\|u\|^h_{\dot{B}_{p,1}^{\frac np+1}}\ge C\big(\|u\|^h_{\dot{B}_{p,1}^{\frac np-1}}\big)^{\beta}, \quad\|\aa\|^h_{\dot{B}_{p,1}^{\frac np}}\ge C\big(\|\aa\|^h_{\dot{B}_{p,1}^{\frac np}}\big)^{\beta}.
\end{align}

Thus, to get the Lyapunov-type inequality, we have to control  $\|( u,\aa)\|^\ell_{\dot{B}_{2,1}^{\frac n2+1}}$ with $\big(\|( u,\aa)\|^\ell_{\dot{B}_{2,1}^{\frac n2-1}}\big)^{\eta}$ for some $\eta>1.$ This process can be obtained  from the interpolation inequality, which implies that we must provide a low order  estimates such as $\|(u,\aa)\|^\ell_{\dot{B}_{2,1}^{-s}}, $  with $-s<\frac n2-1$.
However,  only the incompressible part of stress tensor $(\p\div\tau) $ have dissipation while the whole $\tau$ itself don't. Hence, it's impossible to control $\|(u,\aa)\|^\ell_{\dot{B}_{2,1}^{-s}}$ directly due to  couple terms
$u\cdot \nabla \tau$  and $ F(\tau, \nabla u).$ To overcome this difficulty, we shall control  $\|(u,\tau)\|_{\dot{B}_{2,1}^{-s}}$
instead of $\|(u,\aa)\|^\ell_{\dot{B}_{2,1}^{-s}}$. The price we have to pay is that the stronger condition imposed on $(u_0,\tau_0)\in{\dot{B}_{2,1}^{-s}}(\R^n)$ instead of
$(u_0,\aa_0)\in{\dot{B}_{2,1}^{-s}}(\R^n).$

To do this we
apply   $\dot{\Delta}_j$ to the first two equations in  \eqref{m} and use a standard commutator's process  to get
\begin{eqnarray}\label{R1}
\left\{\begin{aligned}
&\partial_t \dot{\Delta}_ju+ u\cdot\nabla \dot{\Delta}_ju+\ddj\nabla{p}-\Delta \dot{\Delta}_ju-\dot{\Delta}_j\div \tau=[u\cdot\nabla, \dot{\Delta}_j]u,\nonumber\\
&
\partial_t \dot{\Delta}_j\tau+ u\cdot \nabla \dot{\Delta}_j\tau   +\dot{\Delta}_j F(\tau, \nabla u)  -  \dot{\Delta}_jD (u)=[u\cdot\nabla, \dot{\Delta}_j]\tau.
\end{aligned}\right.
\end{eqnarray}

Taking $L^2$ inner product with $\dot{\Delta}_ju, \dot{\Delta}_j\tau$, respectively and using
the following cancellations
$$\int_{\R^n}\dot{\Delta}_j\div \tau\cdot \dot{\Delta}_ju\ dx+\int_{\R^n} \dot{\Delta}_jD (u)\cdot \dot{\Delta}_j\tau \ dx=0,$$
$$
\int_{\R^n} u\cdot\nabla \dot{\Delta}_ju\cdot\dot{\Delta}_ju\ dx=\int_{\R^n} \dot{\Delta}_j\nabla{p}\cdot\dot{\Delta}_ju\ dx=\int_{\R^n} u\cdot\nabla \dot{\Delta}_j\tau\cdot\dot{\Delta}_j\tau\ dx=0,$$
we have
\begin{align*}
\frac{1}{2}\frac{d}{dt}(\|\dot{\Delta}_ju\|_{L^2}^2+\|\dot{\Delta}_j\tau\|_{L^2}^2)
\le& C\|[u\cdot\nabla, \dot{\Delta}_j] u\|_{L^2}\|\dot{\Delta}_j u\|_{L^2}
\nonumber\\
&+C\|\dot{\Delta}_j F(\tau, \nabla u)\|_{L^2}\|\dot{\Delta}_j\tau\|_{L^2}+C\|[u\cdot\nabla, \dot{\Delta}_j] \tau\|_{L^2}\|\dot{\Delta}_j \tau\|_{L^2},
\end{align*}
which implies that
\begin{align}\label{R3}
\frac{d}{dt}\|(\dot{\Delta}_ju,\dot{\Delta}_j\tau)\|_{L^2}
\le& C\big(\|[u\cdot\nabla, \dot{\Delta}_j] u\|_{L^2}
+\|[u\cdot\nabla, \dot{\Delta}_j] \tau\|_{L^2}
+\|\dot{\Delta}_j F(\tau, \nabla u)\|_{L^2}\big).
\end{align}

Integrating  the above inequality from $0$ to $t$,  and multiplying   by $2^{-js}$, we get by summing up about $j\in\Z $ that
\begin{align}\label{R4}
\|(u,\tau)(t,\cdot)\|_{\dot{B}_{2,1}^{-s}}\le&\|(u_0,\tau_0)\|_{\dot{B}_{2,1}^{-s}}+C\int_0^t\| F(\tau, \nabla u)\|_{\dot{B}_{2,1}^{-s}}\ dt'\nonumber\\
&+C\int_0^t(\sum_{j\in\Z}2^{-js}\|[\ddj,u\cdot\nabla] u\|_{L^2} +\sum_{j\in\Z}2^{-js}\|[\ddj,u\cdot\nabla] \tau\|_{L^2})\ dt'.
\end{align}
For any $-\frac np\le s<\frac np$,
from Lemma \ref{daishu} and Lemma \ref{jiaohuanzi}, one has
\begin{align}\label{R5}
\| F(\tau, \nabla u)\|_{\dot{B}_{2,1}^{-s}}\le C\|\nabla u\|_{\dot{B}_{p,1}^{\frac np}}\|\tau\|_{\dot{B}_{2,1}^{-s}}
\le C(\| u\|^\ell_{\dot{B}_{2,1}^{\frac n2+1}}+\| u\|^h_{\dot{B}_{p,1}^{\frac np+1}})\|\tau\|_{\dot{B}_{2,1}^{-s}},
\end{align}
\begin{align}\label{R6}
&\sum_{j\in\Z}2^{-js}\|[\ddj,u\cdot\nabla] u\|_{L^2} +\sum_{j\in\Z}2^{-js}\|[\ddj,u\cdot\nabla] \tau\|_{L^2}\nonumber\\
&\quad\le C\|\nabla u\|_{\dot{B}_{p,1}^{\frac np}}\|u\|_{\dot{B}_{2,1}^{-s}}+\|\nabla u\|_{\dot{B}_{p,1}^{\frac np}}\|\tau\|_{\dot{B}_{2,1}^{-s}}\nonumber\\
&\quad\le C(\| u\|^\ell_{\dot{B}_{2,1}^{\frac n2+1}}+\| u\|^h_{\dot{B}_{p,1}^{\frac np+1}})\|(u,\tau)\|_{\dot{B}_{2,1}^{-s}}.
\end{align}

Plugging the above two estimates into \eqref{R4} implies
\begin{align}\label{A4}
\|(u,\tau)(t,\cdot)\|_{\dot{B}_{2,1}^{-s}}\le&\|(u_0,\tau_0)\|_{\dot{B}_{2,1}^{-s}}\nonumber\\
&+C\int_0^t(\| u\|^\ell_{\dot{B}_{2,1}^{\frac n2+1}}+\| u\|^h_{\dot{B}_{p,1}^{\frac np+1}})\|(u,\tau)\|_{\dot{B}_{2,1}^{-s}}\ dt'.
\end{align}
It is easy to  deduce from
the definition of ${\mathcal X}(t)$ in Theorem \ref{th1} that
\begin{align*}
&\int_0^t(\| u\|^\ell_{\dot{B}_{2,1}^{\frac n2+1}}+\| u\|^h_{\dot{B}_{p,1}^{\frac np+1}})\ dt'\le{\mathcal X}_0.
\end{align*}

Hence,
by the Gronwall inequality, one  can get from \eqref{A4}, for any $-\frac np\le s<\frac np$,  that
\begin{align}\label{A5}
\|(u,\tau)(t,\cdot)\|_{\dot{B}_{2,1}^{-s}}\le C_0
\end{align}
for all $t \ge 0$, where $C_0 > 0$ depends on the norm of $\|(u_0,\tau_0)\|_{\dot{B}_{2,1}^{-s}}$ and ${\mathcal X}_0$.

For any $ s>1-\frac n2,$
it follows from interpolation inequality in Lemma \ref{qianru}  that
\begin{align*}
\|( u,\aa)\|^\ell_{\dot{B}_{2,1}^{\frac n2-1}}
\le& C \big(\|( u,\aa)\|^\ell_{\dot{B}_{2,1}^{-s}}\big)^{\theta_{1}}\big(\|( u,\aa)\|^\ell_{\dot{B}_{2,1}^{\frac n2+1}}\big)^{1-\theta_{1}}\nonumber\\
\le&C\big(\|( u,\tau)\|^\ell_{\dot{B}_{2,1}^{-s}}\big)^{\theta_{1}}\big(\|( u,\aa)\|^\ell_{\dot{B}_{2,1}^{\frac n2+1}}\big)^{1-\theta_{1}},
\quad \theta_1=\frac{4}{n+2s+2}\in(0,1),
\end{align*}
this together with \eqref{A5} implies that
\begin{align}\label{A7}
\|( u,\aa)\|^\ell_{\dot{B}_{2,1}^{\frac n2+1}}\ge  C\big(\|( u,\aa)\|^\ell_{\dot{B}_{2,1}^{\frac n2-1}}\big)^{\frac{1}{1-\theta_{1}}}.
\end{align}
Taking $\beta={\frac{1}{1-\theta_{1}}}$ in \eqref{A812} gives
\begin{align}\label{A8}
\|u\|^h_{\dot{B}_{p,1}^{\frac np+1}}\ge C\big(\|u\|^h_{\dot{B}_{p,1}^{\frac np-1}}\big)^{\frac{1}{1-\theta_{1}}}, \quad\|\aa\|^h_{\dot{B}_{p,1}^{\frac np}}\ge C\big(\|\aa\|^h_{\dot{B}_{p,1}^{\frac np}}\big)^{\frac{1}{1-\theta_{1}}}.
\end{align}

Thus, inserting \eqref{A7} and \eqref{A8} into \eqref{A2} yields
\begin{align*}
\frac{d}{dt}\big(\|(u,\aa)\|^\ell_{\dot{B}_{2,1}^{\frac n2-1}}+\|u\|^h_{\dot B^{\frac  np-1}_{p,1}}+\|\aa\|^h_{\dot B^{\frac  np}_{p,1}}\big)+\bar{c}\big(\|(u,\aa)\|^\ell_{\dot{B}_{2,1}^{\frac n2-1}}+\|u\|^h_{\dot B^{\frac  np-1}_{p,1}}+\|\aa\|^h_{\dot B^{\frac  np}_{p,1}}\big)^{\frac{ n+2s+2}{ n+2s-2}}\le 0.
\end{align*}
Solving this differential inequality directly, we obtain
\begin{align*}
\|(u,\aa)\|^\ell_{\dot{B}_{2,1}^{\frac n2-1}}+\|u\|^h_{\dot B^{\frac  np-1}_{p,1}}+\|\aa\|^h_{\dot B^{\frac  np}_{p,1}}\le&C ({\mathcal X}_0^{-\frac{4}{n+2s-2}}+\frac{4\bar{c}}{n+2s-2}t)^{-\frac{n+2s-2}{4}}\\
\le& C(1+t)^{-\frac{n+2s-2}{4}}.
\end{align*}
Moreover, from Lemma \ref{qianru}, we further get
\begin{align}\label{A10}
\|(u,\aa)\|_{\dot{B}_{p,1}^{\frac np-1}}\le C(\|(u,\aa)\|^\ell_{\dot{B}_{2,1}^{\frac n2-1}}+\|u\|^h_{\dot{B}_{p,1}^{\frac np-1}}+\|\aa\|^h_{\dot{B}_{p,1}^{\frac np}})\le C(1+t)^{-\frac{n+2s-2}{4}}.
\end{align}

For any $\frac np-\frac n2-s<\gamma<\frac np-1,$ by the interpolation inequality we have
\begin{align*}
\|(u,\aa)\|^\ell_{\dot{B}_{p,1}^{\gamma}}
\le& C\|(u,\aa)\|^\ell_{\dot{B}_{2,1}^{\gamma+\frac n2-\frac np}}\\
\le&C\big(\|(u,\aa)\|^\ell_{\dot{B}_{2,1}^{-s}}\big)^{\theta_{2}} \big(\|(u,\aa)\|^\ell_{\dot{B}_{2,1}^{\frac n2-1}}\big)^{1-\theta_{2}},\quad \theta_{2}=\frac{\frac np -1-\gamma}{\frac n2-1+s}\in (0,1),
\end{align*}
which combines \eqref{A5} with \eqref{A10} gives
\begin{align}\label{A10234}
\|(u,\aa)\|^\ell_{\dot{B}_{p,1}^{\gamma}}
\le C(1+t)^{-\frac{(\frac n2+s-1)\theta_{2}}{2}}=C(1+t)^{-\frac{n}{2}(\frac 12-\frac 1p)-\frac{s+\gamma}{2}}.
\end{align}
In  the light of
$\frac np-\frac n2-s<\gamma<\frac np-1,$
 we see that
$$\|(u^h,\aa^h)\|_{\dot{B}_{p,1}^{\gamma}}\le C(\|u\|^h_{\dot{B}_{p,1}^{\frac np-1}}+\|\aa\|^h_{\dot{B}_{p,1}^{\frac np}})\le C(1+t)^{-\frac{n+2s-2}{4}},
$$
from which and \eqref{A10234} gives
\begin{align*}
\|(u,\aa)\|_{\dot{B}_{p,1}^{\gamma}}
\le&C(\|(u,\aa)\|^\ell_{\dot{B}_{p,1}^{\gamma}}+\|(u,\aa)\|^h_{\dot{B}_{p,1}^{\gamma}})\nonumber\\
\le& C(1+t)^{-\frac{n}{2}(\frac 12-\frac 1p)-\frac{s+\gamma}{2}}+C(1+t)^{-\frac{n+2s-2}{4}}\nonumber\\
\le& C(1+t)^{-\frac{n}{2}(\frac 12-\frac 1p)-\frac{s+\gamma}{2}}.
\end{align*}

Thanks to the embedding relation
$\dot{B}^{0}_{p,1}(\R^n)\hookrightarrow L^p(\R^n)$, one infer that
\begin{align*}
\|\Lambda^{\gamma} (u,\aa)\|_{L^p}
\le& C(1+t)^{-\frac{n}{2}(\frac 12-\frac 1p)-\frac{s+\gamma}{2}}.
\end{align*}
For any
$p\leq q\leq\infty$  and $\frac nq-\frac np-s<\alpha \leq\frac nq-1$, by the Gagliardo-Nirenberg type interpolation inequality,
which can be found in the  Chap. 2 of \cite{bcd},
taking
$$
k\theta_{3}+m(1-\theta_{3})=\alpha+n\Bigl(\frac1p-\frac1q\Bigr), \quad m=\frac np-1,
$$
we get
\begin{align*}
\|\Lambda^{\alpha}(u,\aa)\|_{L^{q}} \le &C\|\Lambda^{m}(u,\aa)\|_{L^p}^{1-\theta_{3}}\|\Lambda^{k}(u,\aa)\|^{\theta_{3}}_{L^p}
\nonumber\\ \le & C\Big\{(1+t)^{-\frac n2(\frac 12-\frac 1p)-\frac{m+s}{2}}\Big\}^{1-\theta_{3}}
\Big\{(1+t)^{-\frac n2(\frac 12-\frac 1p)-\frac{k+s}{2}}\Big\}^{\theta_{3}}
\nonumber\\=& C(1+t)^{-\frac n2(\frac 12-\frac 1q)-\frac {\alpha+s}{2}}.
\end{align*}

Consequently, we have completed the proof of our theorem. \quad\quad$\Box$

\bigskip
\noindent \textbf{Acknowledgement.} This work is supported by NSFC under grant number 11601533.



\begin{thebibliography}{99}



\bibitem{bcd}
H.~Bahouri, J.Y. Chemin, R.~Danchin,
\newblock { {F}ourier {A}nalysis and {N}onlinear {P}artial {D}ifferential
  {E}quations}.
\newblock Grundlehren Math. Wiss. , vol. {\textbf{343}}, Springer-Verlag,
  Berlin, Heidelberg, 2011.
\newblock $\,$

\bibitem{bird}
R.B. Bird, C.F. Curtiss, R.C. Armstrong, O. Hassager,
\newblock Dynamics of polymeric liquids.
\newblock {Fluid Mechanics, vol. 1,  2nd edn
Wiley, New York, 1987.}
\newblock $\,$


\bibitem{chemin}
J.Y. Chemin, N. Masmoudi,
\newblock About lifespan of regular solutions of equations related to viscoelastic fluids.
\newblock {\it SIAM J. Math. Anal.}, {\bf 33}, 84--112, 2001.
\newblock $\,$

\bibitem{chenqionglei}
Q. Chen, X. Hao,
\newblock Global well-posedness in the critical Besov spaces for the incompressible Oldroyd-B model without damping mechanism.
\newblock {\it ArXiv:1810.06171.}
\newblock $\,$

\bibitem{miaochangxing}
Q. Chen, C. Miao,
\newblock Global well-posedness of viscoelastic fluids of Oldroyd type in Besov spaces.
\newblock {\it Nonlinear Anal.}, {\bf 68}, 1928-1939, 2008.
\newblock $\,$





\bibitem{constanin}
P. Constantin, M. Kliegl,
\newblock Note on global regularity for two-dimensional Oldroyd-B fluids with diffusive stress.
\newblock {\it Arch. Ration. Mech. Anal.}, {\bf 206}, 725--740, 2012.
\newblock $\,$


\bibitem{EL}
T. M. Elgindi,  J. Liu,
{Global wellposedness to the generalized Oldroyd type models in $\mathbb{R}^3$.}
\newblock {\it  J. Differential Equations}, {\bf 259},  1958--1966, 2015.
\newblock $\,$


\bibitem{ER}
T.M. Elgindi,  F. Rousset,
{Global regularity for some Oldroyd-B type models}.
\newblock {\it  Comm. Pure Appl. Math.}, {\bf 68},  2005--2021, 2015.
\newblock $\,$


\bibitem{FZ}
D. Fang, R. Zi,
{Global solutions to the Oldroyd-B model with a class of large initial data}.
\newblock {\it  SIAM J. Math. Anal.}, {\bf 48},  1054--1084, 2016.
\newblock $\,$


\bibitem{GS}
C. Guillop\'{e}, J.C. Saut,
{Existence results for the flow of viscoelastic fluids with a differential constitutive law}.
\newblock {\it  Nonlinear Anal.}, {\bf 15},  849--869, 1990.
\newblock $\,$





\bibitem{GS2}
C. Guillop\'{e}, J.C. Saut,
{Global existence and one-dimensional nonlinear stability of shearing motions of viscoelastic fluids of Oldroyd type}.
\newblock {\it  RAIRO Mod\'{e}l. Math. Anal. Num\'{e}r.}, {\bf 24},  369--401, 1990.
\newblock $\,$

\bibitem{guoyan}
Y. Guo, Y. Wang,
{Decay of dissipative equations and negative Sobolev spaces}.
\newblock {\it  Comm. Partial Differential Equations}, {\bf 37},  2165--2208, 2012.
\newblock $\,$





\bibitem{LMZ}
Z. Lei, N. Masmoudi,  Y. Zhou,
{Remarks on the blowup criteria for Oldroyd models}.
\newblock {\it  J. Differential Equations}, {\bf 248}, 328--341, 2010.
\newblock $\,$






\bibitem{lin2012}
F. Lin,
\newblock  Some analytical issues for elastic complex fluids.
\newblock {\it Comm. Pure Appl. Math.}, {\bf 65}, 893--919, 2012.
\newblock $\,$



\bibitem{LM}
P.L. Lions,  N. Masmoudi,
{Global solutions for some Oldroyd models of non-Newtonian flows}.
\newblock {\it  Chinese Ann. Math. Ser. B}, {\bf 21}, 131--146, 2000.
\newblock $\,$

\bibitem{Oldroyd}
J. Oldroyd,
{Non-Newtonian effects in steady motion of some idealized elastico-viscous
liquids}.
\newblock {\it  Proc. Roy. Soc. Edinburgh Sect. A}, {\bf 245}, 278--297, 1958.
\newblock $\,$


\bibitem{xinzhouping}
Z. Xin, J. Xu,
{Optimal decay for the compressible Navier-Stokes equations without additional smallness assumptions}.
\newblock {\it  ArXiv:1812.11714v1}.
\newblock $\,$



\bibitem{zhaixiaoping}
X. Zhai,
\newblock Global  solutions to the  $n$-dimensional incompressible Oldroyd-B model without damping mechanism.
\newblock {\it  ArXiv:1810.08048}.
\newblock $\,$

\bibitem{zhuyi}
Y. Zhu,
\newblock Global small solutions of 3D incompressible Oldroyd-B model without damping mechanism.
\newblock {\it  J. Funct. Anal.}, {\bf 274}, 2039--2060, 2018.
\newblock $\,$


\bibitem{zuiruizhao}
R. Zi, D. Fang, T. Zhang,
\newblock Global solution to the incompressible
Oldroyd-B model in the critical $L^p$
framework: the case of the non-small
coupling parameter.
\newblock {\it  Arch. Rational Mech. Anal.}, {\bf 213},  651--687, 2014.
\newblock $\,$













\end{thebibliography}
\end{document}